\documentclass[11pt]{article}
\usepackage{stmaryrd}
\usepackage{}
\usepackage{bbm}

\usepackage{mathrsfs}
\usepackage{amsfonts}

\renewcommand{\baselinestretch}{1.15}

\renewcommand{\baselinestretch}{1.15}

\usepackage{amsmath}
\usepackage{amsxtra}
\usepackage{amstext}
\usepackage{amssymb}
\usepackage{latexsym}
\usepackage{array}
\usepackage{dsfont}
\usepackage{graphicx}
\usepackage{cite}
\usepackage{color}
\usepackage{subfigure}

\newtheorem{Th}{Theorem}[section]
\newtheorem{aTh}{Auxiliary Theorem}[section]
\newtheorem{Le}{Lemma}[section]
\newtheorem{Co}{Corollary}[section]
\newtheorem{De}{Definition}[section]

\newtheorem{Pro}{Proposition}[section]
\newtheorem{Rem}{Remark}[section]
\newtheorem{Ex}{Example}[section]
 \newtheorem{Con}{Conjecture}[section]
 \newtheorem{Concl}{Conclusion}[section]

\newcommand{\bth}{\begin{Th}}
\newcommand{\eeth}{\end{Th}}
\newcommand{\bath}{\begin{aTh}}
\newcommand{\eeath}{\end{aTh}}
\newcommand{\ble}{\begin{Le}}
\newcommand{\ele}{\end{Le}}
\newcommand{\bco}{\begin{Co}}
\newcommand{\eco}{\end{Co}}
\newcommand{\bde}{\begin{De}}
\newcommand{\ede}{\end{De}}
\newcommand{\bpr}{\begin{Pro}}
\newcommand{\epr}{\end{Pro}}
\newcommand{\bre}{\begin{Rem}}
\newcommand{\ere}{\end{Rem}}
\newcommand{\bex}{\begin{Ex}}
\newcommand{\eex}{\end{Ex}}
\newcommand{\beq}{\begin{equation}}
\newcommand{\eeq}{\end{equation}}
\newcommand{\beqn}{\begin{eqnarray}}
\newcommand{\eeqn}{\end{eqnarray}}
\newcommand{\be}{\begin{eqnarray*}}
\newcommand{\ee}{\end{eqnarray*}}
\newcommand{\m}[1]{\mbox{\bf #1}}
\newcommand{\p}{\m{Proof:$\;$}}

\newcommand{\ef}{\begin{flushright}\vspace{-8.2mm} $\!\blacksquare$  \end{flushright}}
\newcommand{\bconcl}{\begin{Concl}}
\newcommand{\econcl}{\end{Concl}}
\newcommand{\bcon}{\begin{Con}}
\newcommand{\econ}{\end{Con}}

\def\mysection#1{\setcounter{equation}{0}
{\renewcommand{\thesection}{\arabic{section}}\section{#1}}
\renewcommand{\thesection}{\arabic{section}}}

\begin{document}
\title{
{\huge \bf Hilbert Boundary Value Problems for Hyper Monogenic Functions on The Hyperplane  }\!
\thanks{\setlength{\baselineskip}{12pt}
This work was
 partially supported by the National Natural Science Foundation of China (Grant No.
 11171260, 111713.35), the Science and Technology
 Development Fund, Macao SAR (Grant No. 0006/2019/A1,  0123/2018/A3), Natural Science Foundation of Fujian Province (2020J01322)
   DYSP and AMEP of Linyi University.}}
\author{\renewcommand{\baselinestretch}{1}
 Pei Dang$^{\,a}$,\,\,\, Jinyuan Du$^{\,b,\,c,\,}$\footnote{Corresponding author. Email: jydu@whu.edu.cn},\,\,\,\,\, Tao Qian$^{\,d}$\\
\small a. Faculty of Information  Technology, Macau University of Science and Technology, Macao\\
\small b. Department of Mathematics, Wuhan University, Wuhan $430072$, P. R. China\\
\small c. School of Science, Linyi University, Linyi, Shandong 276000, P. R. China\\
\small d. {Faculty of Information  Technology, Macau University of Science and Technology, Macao}
\vspace{-9mm}}
\date{}
\maketitle


\begin{center}
\begin{minipage}{131mm}
\mbox{}\hspace{4mm}
{\small\bf Abstract.} {\small {This paper systematically studies Hilbert boundary value problems for hyper monogenic functions on the hyperplane for the solutions being of any integer orders at the infinity, where the negative order cases are new even when restricted to the complex plane context. The explicit solution formulas are given and the solvability conditions are specified. The results are proved through using the Clifford symmetric extension method to reduce Hilbert boundary value problems to Riemann boundary value problems.}
}

\mbox{}\hspace{4mm}{\small\bf Keywords}\,\,\, {\small Real part and imaginary part; Order at the infinite point; Hilbert boundary value problem; Symmetric extension.}

\mbox{}\hspace{4mm}{\small\bf 1991 MR Subject Classification}\, {\small  30E20,
 30E25, 30G35, 30G30,   31B10.}
\end{minipage}
\end{center}

\mysection{Introduction}

 {As is well known, boundary value problems (BVPs) for analytic functions in the classical
complex analysis form an important branch of mathematical analysis, and, due to its theoretical elegance and ample applications in physics and other
subjects such as elastic theory, fracture mechanics, hydromechanics, etc., have been thoroughly studied over a long time.}  The theory of boundary
value problems for analytic functions has been investigated
systematically in {the literature}\cite{Lu, Mu, Ga}, and  applications have been showed roundly in {the monographs} \cite{lu2,mu2}.

It has been proved that the function theory over Clifford algebra is
an appropriate setting to generalize many aspects of {the function} theory
of one complex variable to higher dimensions
\cite{BDS-RNM-82,DelangheSS-92,GurlebeckS-90,GurlebeckS-97}.
 It is natural that mathematicians hope to develop theories
on boundary value problems for hyper monogenic functions, also called regular functions simply,  in the hypercomplex
analysis analogous to those for analytic functions in the classical
complex analysis. In fact, some results on boundary value problems
for analytic functions in the classical complex analysis have been
generalized to {regular functions} in Clifford analysis ($e.g.$, some articles listed in
 \cite{dd}).

In \cite{ld} and \cite{ddq17}, we  discussed, respectively, Riemann boundary value problems on {closed  smooth  surfaces} and on  the {hyperplanes} in detail. In this paper, we consider  Hilbert boundary value problems {(Hilbert BVPs) on the half-hyperplanes}, such as the Poincar\'{e} upper half space. {Hilbert BVPs} are  very important in pure mathematics and engineering  practice. The discussion is rather difficult and complex. Even in the classical complex analysis,  Hilbert boundary value problems on the real axis are also  not completely discussed \cite{Lu, Mu}. The researchers are restricted to the case of the bounded growth at the infinity, {while for general growth orders at the infinity there have some obstacles}.
In \cite{Xu}, Z. Y. Xu  and C. P. Zhou tried to solve Hilbert boundary value problems on the hyperplane.
They generalized the classical results for Hilbert boundary value problems on the real axis to the hypercomplex analysis setting under the condition that the solutions are bounded  at the infinity.
In \cite{GD}, Y. F. Gong and J. Y. Du {continued} to discuss  Hilbert BVPs  under the condition that the solutions have  finite {non-negative orders}  at the infinity  by the symmetric {extension method}. {The} arguments of \cite{Xu,GD} {contain} some mistakes that will be pointed out below. {To the authors' knowledge, Hilbert BVPs with negative order at infinity have not been discussed in any paper, because the corresponding Riemann BVPs were not discussed before \cite{ddq17}, so Hilbert BVPs with negative order at infinity are always an open problem.}

In this paper, we will systematically discuss the Hilbert boundary value problems
on the hyperplane for regular functions in Clifford analysis, including the {cases of negative orders at the infinity}.
 We will use the so-called {symmetric extension method
to solve the Hilbert BVPs}.
The paper is organized as follows.
In \S 2 we review some of the necessary preparations in Clifford analysis.
In \S 3, we formulate  the  Hilbert boundary value problems in the Poincar\'{e} upper half space.
 In \S 4, we introduce the symmetric extension of a regular function in the Poincar\'{e} upper half space and discuss its regularity. Then the Hilbert BVPs  are sloved by converting them equivalently {to the} Riemann BVPs discussed in \cite{ddq17}.
 In \S 5, {the solution formulas and the conditions of the Schwarz BVPs and the Hilbert BVPs for all orders at the infinity of the solutions are  obtained}.
  These results  {extend} both the classical ones in complex analysis
\cite{Lu, Mu, Ga} and those in the Clifford analysis setting \cite{Xu, GD}.
\vspace{2mm}

\mysection {Hypercomplex functions }

 We begin by {recalling} the necessary preliminary knowledge {in} Clifford algebra and Clifford analysis \cite{DelangheSS-92, BDS-RNM-82}, which are used throughout this paper.
\vspace{2.5mm}

\subsection{Clifford analysis}
\mbox{}\hspace{6mm}
Let $C(V_{n})$  be a  $2^n$-dimensional real linear space. To expediently introduce the product on it, we write
its basis by
$\left\{ e_A, A=(h_1,\cdots,
h_r)\!\in\!\mbox{$\cal P$}\!N, 1\leq h_1<\cdots<h_r \leq n
\right\}$,  where $N$ stands for the set
$\{1,\cdots, n \}$ and $\mbox{$\cal P$}\!N$ denotes the
family of all the ordered  subsets of $N$ in the above fixed
way.
 Sometimes,
 $e_{\emptyset}$ is written as $e_0$ and
 $e_A$ as $e_{h_1\cdots h_r}$ for $A=\{h_1,\cdots,h_r\}\in\mbox{$\cal P$}\!N$.
The product on $C(V_{n})$ is defined by
\beq \label{2.1}
\left\{
\begin{array}{ll}
e_Ae_B=(-1)^{\#(A\cap B)}(-1)^{P(A,B)}e_{A\mathbf{\Delta}
B},
&\mbox{\it if}\,\, \,\,A, B\in\mbox{$\cal P$}\!N, \\[4.5mm]
\lambda\mu =\!\sum\limits_{A\in\mbox{\small$ \cal P$}\!
\mbox{\small$ N$}}\sum\limits_{B\in\mbox{\small$ \cal P$}
\!\mbox{\small$ N$}}\lambda_{A}\mu_{B}e_Ae_B, & \mbox{\it
if}\,\,\,\, \lambda\!=\!\! \sum\limits_{A\in\mbox{\small$ \cal
P$}\!\mbox{\small$ N$}}\lambda_Ae_A,\,\,\,
\mu\!=\!\!\sum\limits_{A\in\mbox{\small$ \cal P$}\!\mbox{\small$
N$}} \mu_{A}e_A,
\end{array}
\right. \eeq
where the notation
$\#(A)$ denotes  the number of the elements in $A$ and
 $P(A,B)=\sum\limits_{j\in B}P(A,j)$ with $P(A,j)=\#\{i, i\in A, i>j\}$,
the symmetric difference set $A\mathbf{\Delta}B$ is also an ordered one in the above way, and
$\lambda_{A}\in\mathbb{R}$ is the coefficient of the
$e_A$--component of the Clifford number $\lambda$. It follows at
once from the multiplication rule $(\ref{2.1})$ that $e_0$ is the
identity element written {as} $1$ and in particular,
\beq \label{2.2}
\left\{
 \begin{array}{ll}
 e_i^2=-1,&\mbox{\it if}\,\,\,\, \,\,i=1, \cdots, n,\\[1mm]
  e_ie_j=-e_je_i,&\mbox{\it if}\,\,\,\, \,\, 1\leq i<j\leq n,\\[1mm]
 e_{h_1}e_{h_2}\cdots e_{h_r}=e_{h_1h_2\cdots h_r},\hspace{3mm}&\mbox{\it if}\,\, \,\,\,\,1\leq h_1<h_2<\cdots
<h_r\leq n.
\end{array}
\right. \eeq

It is clear that $C({V}_{n})$ is a real linear and associative,
non-commutative algebra by algebraically spanning the linear
subspace $V_{n}=\mbox{\rm span}\{e_1,e_2,\cdots,e_n\}$.
 It is called the Clifford algebra over
$V_{n}$. The elements $\lambda=\lambda_0+\lambda_1e_1+\cdots+\lambda_n e_n$ for $\lambda_0,\cdots,\lambda_n\in\mathbb{R}$ are called paravectors.

We frequently use the following defined involution:
\beq
\label{2.3}
\left\{
\begin{array}{ll}
\overline{e_A}=(-1)^{\frac{\#(A)(\#(A)+1)}{2}}e_A, &\mbox{\it if}\,\, \,\,\,\,A\in\mbox{$\cal P$}\!N, \\[4mm]
\overline{\lambda}=\sum\limits_{A\in\mbox{\small$ \cal
P$}\!\mbox{\small$ N$}}\lambda_{A} \overline{e_{A}},
\hspace{20mm} & \mbox{\it if}\,\,\,\,\,\,
\lambda=\sum\limits_{A\in\mbox{\small $\cal P$}\!\mbox{\small$
N$}}\lambda_{A} e_{A}.
\end{array}
\right.
\eeq
In the sequel,
$\lambda_A$ is also written as $[\lambda]_A$. In particular, the
coefficient $\lambda_{\emptyset}$ is denoted by $\lambda_0$ or
$[\lambda]_0$, which is called the scalar part of the
 Clifford number $\lambda$.
An inner product $(\cdot,\cdot)$ on $C\left(V_{n}\right)$ is
defined by putting for any $\lambda$ and $\mu$ in
$C\left(V_{n}\right)$
\beq \label{2.4}
\big(\lambda,\mu\big)=\big[\lambda\overline{\mu}\big]_0=\sum\limits_A\lambda_A\mu_A,
 \eeq where
$\lambda=\sum\limits_{A}\lambda_Ae_A$, $\mu=\sum\limits_A\mu_Ae_A$
and the symbol $\sum\limits_A$ is an abbreviation of
$\sum\limits_{A\in\mbox{\small $\cal
P$}\!\mbox{\small$N$}}$. \vspace{1mm}Thus, the corresponding norm on $C(V_{n})$ reads,
 \beq \label{2.5}
 \big|\lambda\big|=\sqrt{(\lambda,\lambda)}=
\left[\sum\limits_A\lambda_A^2\right]^{\frac{1}{2}}.
 \eeq

In such way, $C(V_{n})$ is a real Hilbert space and at the same
time it is a Banach algebra with the equivalent norm
\beq
\label{2.6}
\big|\lambda\big|_0=2^{\frac{n}{2}}\big|\lambda\big|,
\eeq
 that is
\beq \label{2.7}
\big|\lambda\mu\big|_{0}\leq\big|\lambda\big|_{0}\big|\mu\big|_{0},\,\,\,\,\big|\lambda\mu\big|\leq 2^{\frac{n+1}{2}}\big|\lambda\big|\big|\mu\big|. \eeq
In  particular, if $\lambda$ is a paravector and $\mu\in C(V_n)$, then \cite{ BDS-RNM-82}
\beq \label{2.8}
\big|\lambda\mu\big|=\big|\mu\lambda\big|=\big|\lambda\big|\big|\mu\big|.
\eeq

Let $\Omega$ be a  non-empty subset of $\mathbb{R}^{n+1}$.
Hypercomplex functions $f$ defined in $\Omega$ and with values in
$C\left(V_{n}\right)$ will be considered, $i.e.$,
$
f: \Omega\longrightarrow C\!\left(V_{n}\right).
$
They are of the form
\beq\label{2.9}
f(w)=\sum\limits_{A}f_A(w)e_A,\,\,\,w=\left(w_0,w_1,\cdots,w_n\right)
\in \Omega\subset\mathbb{R}^{n+1},
 \eeq
 where the $f_A(w)$ is the $e_A$--component of $f(w)$.
Obviously, the $f_A$'s are real--valued functions in $\Omega$, which
are called the $e_A$--component functions of $f$. Whenever a
property such as differentiability and continuity is ascribed
to $f$,  it is clear that in fact all the component functions $f_A$
possess the cited property. \vspace{1mm}So the meaning $f\in C^{(r)}\!\left(\Omega,
C\left(V_{n}\right)\right)$ is very clear.

Obviously, $C(V_{n-1})$ is a subalgebra of $C(V_{n})$ where $V_{n-1}=\mbox{\rm span}\{e_1,e_2,\cdots,e_{n-1}\}$.  Then, $\lambda\in C(V_{n})$ has the
unique  decomposition \cite{Xu,GD}
\beq \label{2.10}
\lambda=x+ e_n\,y^l\,\,\,\,\, \mbox{\it where}\,\,\,\,\,x,y^l\in C(V_{n-1}),
\eeq
$i.e.$,
\beq \label{2.11}
C\left(V_{n}\right)=C\left(V_{n-1}\right)\oplus e_n\,C\left(V_{n-1}\right).
\eeq
We define
\beq \label{2.12}
\mathrm{Re} (\lambda)=x,\hspace{8mm} \mathrm{Im}^l\,(\lambda)=y^l.
\eeq

\bre\label{re 2.1}
{\rm Similarly, $\lambda\in C(V_{n})$ has also the
unique  decomposition
\beq \label{2.13}
\lambda=x+ y^r\,e_n\,\,\, \,\,\mathrm{\it where}\,\,\,\,\,x,y^r\in C(V_{n-1}).
\eeq
We  also define
\beq \label{2.14}
\mathrm{Re} (\lambda)=x,\hspace{8mm} \mathrm{Im}^r\,(\lambda)=y^r.
\eeq
For clarity, we call, respectively, $y^l$ and $y^r$ in (\ref{2.10}) and (\ref{2.13}), the left and right imaginary part of
$\lambda$. From (\ref{2.10}) and (\ref{2.13}), we have
\beq\label{2.15}
e_n\,y^l=y^r\,e_n,\,\,\,\,\,\,i.e.,\,\,\,\,\,\, e_n\,\mathrm{Im}^l (\lambda)=
\mathrm{Im}^r (\lambda)\,e_n.
\eeq

} \ere

 It is clear that the decompositions $(\ref{2.10})$ and {$(\ref{2.13})$} are  generalizations of the representation of the classical complex numbers.
In other words,
$(\ref{2.12})$ and $(\ref{2.14})$ are the generalization of operators  $\mathrm{Re}$ and $\mathrm{Im}$ acting on the complex
numbers. From (\ref{2.11}) and (\ref{2.13}) we obviously have
\beq\label{2.16}
\big|\lambda\big|^2=\big|x\big|^2+\big|y^r\big|^2,
\,\,\,\,\,\,\,
\big|\lambda\big|^2=\big|x\big|^2+\big|y^l\big|^2.
\eeq

For a hypercomplex function $f$ given by (\ref{2.9}), we call, respectively,
\beq
\label{2.17}
\big(\mathrm{Re}f\big)(w)=\mathrm{Re}\big(f(w)\big),\,\,\,\,
\,\,\,\,\,\big(\mathrm{Im}^r f\big)(w)=\mathrm{Im}^r\big(f(w)\big),\,\,\,\,\,
\big(\mathrm{Im}^l f\big)(w)=\mathrm{Im}^l\big(f(w)\big),
\eeq
the real part of $f$, the right and left imaginary part of $f$.
In particular,
if
\beq \label{2.18}
\big(\mathrm{Im} f\big)^r =\big(\mathrm{Im}^l f\big)=0,\,\,\,w=\left(w_0,w_1,\cdots,w_n\right)
\in \Omega\subset\mathbb{R}^{n+1},
\eeq
$i.e.$, $
f\!: \Omega\longrightarrow C\!\left(V_{n-1}\right)
$, then we say it  be a $C(V_{n-1})$-valued function, briefly, a para real-valued function which mimics the case of the real-valued function in the classical complex analysis.

Sometimes, for clarity, we write
\beq\label{2.19}
\lambda=\sum\limits_{A\in\mathcal {P} \{1,\cdots,n-1\}}\lambda_{A} e_{A}\,\,\,\,\,\mathit{when}\,\,\,\,\,\lambda\in C\big( {V}_{n-1}\big),
\eeq
where $\mathcal {P} \{1,\cdots,n-1\}$ is to denote the family of all ordered subsets of $\{1,\cdots,n-1\}$  in the similar
 way used in $\mathcal {P} N$. Thus,
a hypercomplex function $f$ given in (\ref{2.9}) may be re-written as

\beq\label{2.20}
\begin{array}{lll}
f(w)&=&\displaystyle\sum_{A\in\mathcal {P} \{1,\cdots,n-1\}}u_{A}(w)\,e_{A}+\sum_{A\in\mathcal {P} \{1,\cdots,n-1\}}v^{l}_{A}(w)e_{A}e_{n}\\[8mm]
&=&\displaystyle\sum_{A\in\mathcal {P} \{1,\cdots,n-1\}}u_{A}(w)e_{A}+
\sum_{A\in\mathcal {P} \{1,\cdots,n-1\}}v^{r}_{A}(w)e_{n}\,e_{A},
\end{array}
\eeq
where $u_A$, $v^{l}_A$ and $v^{r}_A$ are the real-valued functions,\vspace{1mm} which {are called} respectively the component functions on $C\big( V_{n-1}\big)$, $e_n\,C\big(V_{n-1}\big)$ and
$C\big(V_{n-1}\big)e_n$. Obviously,
\beq\label{2.21}
\big(\mathrm{Re}f\big)(w)\triangleq u(w)=\sum_{A\in\mathcal {P} \{1,\cdots,n-1\}}u_{A}(w)\,e_{A},\,\,\,x\in\Omega,
\eeq
and
\beq\label{2.22}
\begin{array}{ll}
\displaystyle\big(\mathrm{Im}f\big)^l (w)&\!\!\!\triangleq v^l(w)=\sum\limits_{A\in\mathcal {P} \{1,\cdots,n-1\}}v^l_{A}(w)\,e_{A},\,\,\,\,\,w\in\Omega,\\[6mm]
\displaystyle\big(\mathrm{Im}f\big)^{r} (w)&\!\!\!\triangleq v^r(w)=\sum\limits_{A\in\mathcal {P} \{1,\cdots,n-1\}}v^r_{A}(w)\,e_{A},\,\,\,\,\,w\in\Omega.
\end{array}
\eeq
Obviously,
\beq\label{2.23}
v^r_A=(-1)^{\#(A)}v^l_A\,\,\,\,\,\,\,\mathrm{\it and}\,\,\,\,\,\,\,e_n\,\big(\mathrm{Im}f\big)^l (w)=\big(\mathrm{Im}f\big)^r (w)\,e_n.
\eeq

When $w=(w_0,w_1,\cdots,w_n)\in\mathbb{R}^{n+1}$, we
  introduce the mapping
 \beq \label{2.24}
 \mathrm{capital}\!:\,w
\longmapsto W=\displaystyle\sum\limits_{i=0}^{n}w_ie_i,
\eeq
 which is {a} proper isomorphism between
$\mathbb{R}^{n+1}$ and the linear subspace $\mbox{\rm span}\{e_0,e_1,\cdots,e_n\}$ of
$C\!\left(V_{n}\right)$. In the sequel, we simply treat
 the capital $W$ as $w$. \vspace{1mm}This is Vahlen's choice \cite{Vahlen, Ahlfors}.

 Thus, we have
  \beq \label{2.25}
\mathrm{Re}\,(w)=w_0+w_1e_{1}+\cdots+w_{n-1}e_{n-1}
\,\,\,\,\,\,\,\,\mbox{\it while}\,\,\,\,\,\,\,\,w=(w_0,w_1,\cdots,w_n)\in\mathbb{R}^{n+1}.
\eeq
and
 \beq \label{2.26}
\mathrm{Im}^{l}(w)=\mathrm{Im}^{r}(w)=w_n\triangleq \mathrm{Im}(w)
\,\,\,\,\,\mbox{\it while}\,\,\,\,\,\,w=(w_0,w_1,\cdots,w_n)\in\mathbb{R}^{n+1}.
\eeq

In the following, if
\beq\label{2.27}
\mathrm{Im}^{l}(\lambda)=\mathrm{Im}^{r}(\lambda),
\eeq
we will write both of them, without confusion, just as
$\mathrm{Im}(\lambda)$.

We define
\beq \label{2.28}
\mathbb{R}^{n+1}_{+}=\big\{ w,\, w\in\mathbb{R}^{n+1},\, \mathrm{Im}(w)>0\big\},\,\,\,\,\,\,\,\,
\mathbb{R}^{n+1}_{-}=\big\{ w,\,w\in\mathbb{R}^{n+1},\, \mathrm{Im}(w)<0\big\},
\eeq
which are called, respectively,  the Poincar\'{e} upper half space and  the Poincar\'{e} lower half space,
 while the hyperplane
 \beq\label{2.29}
 \mathbb{R}^{n+1}_{0}=\big\{ w,\, \mathrm{Im}(w)=0\big\}
 \eeq
 is called the parareal plane in $\mathbb{R}^{n+1}$.

    It is noted that the paravectors given in $\mathbb{R}^{n+1}_{0}$ play a treble role as elements of $\mathbb{R}^{n}$ and $\mathbb{R}^{n+1}_{0}$ as well as  $C(V_{n-1})\subset C( V_n)$.
\vspace{1.5mm}

\subsection{Regular functions}
\mbox{}\hspace{6mm}
Let $\Omega$ be a domain of $\mathbb{R}^{n+1}$. Introduce the following Dirac operator
\beq \label{2.30}
D=\sum\limits^{n}_{k=0}e_{k}\displaystyle\frac{\partial} {\partial
x_{k}}:\,\,C^{(r)}\!\left(\Omega,
C\left(V_{n}\right)\right)\longrightarrow
C^{(r-1)}\!\left(\Omega, C\left(V_{n}\right)\right),
\eeq
whose actions on functions from the left and from the right are
governed by the rules
\beq \label{2.31}
D[f]=\sum\limits_{k=0}^{n}\sum\limits_{A}e_{k}e_A
\displaystyle\frac{\partial f_A} {\partial x_{k}}, \,\,\,\,
[f]D=\sum\limits_{k=0}^{n}\sum\limits_{A}e_Ae_{k}
\displaystyle\frac{\partial f_A}{\partial x_{k}}.
\eeq

\bde \label{de 2.1}
We say that a function $f\in C^{(r)}\!\left(\Omega,
C\left(V_{n}\right)\right) \,\,(r\geq1)$ is left
$($right$)$ regular or monogenic in $\Omega$ if $D[f]\!=\!0$ $([f]D\!=\!0)$ in
$\Omega$. $f$ is said to be biregular in $\Omega$  if
it is both left  and right regular. \ede

\bre \label{re 2.2}
 {\rm  Generally speaking, if
\beq \label{2.32}
f(w)=\sum\limits_{j=1}^{m}f_j(w)\lambda_j\,\,\,\,\Big(f_j\in C(\Omega,\mathbb{R})\,\,\mbox{\it and}\,\,
\lambda_j\!\in\! C(V_{n})\Big)
\eeq
 then
\beq \label{2.33}
D[f]=\displaystyle\sum\limits_{k=0}^{n}\sum\limits_{j=1}^{m}e_{k}\lambda_j
\frac{\partial f_j} {\partial w_{k}},\,\,\,\,\,\,\,
[f]D=\displaystyle\sum\limits_{k=0}^{n}\sum\limits_{j=1}^{m}\lambda_j e_{k}
\frac{\partial f_j} {\partial w_{k}}.
\eeq
 } \ere

\bex\label{ex 2.1}{ Let
\beq \label{2.34}
E(w)=\displaystyle\frac{\overline{w}}{|w|^{n+1}}=\frac{w^{-1}}{|w|^{n-1}}, \hspace{4mm}w\in\mathbb{ R}^{n+1}\setminus\{0\}.
\eeq
Then $E$ is biregular \mbox{\rm \cite{BDS-RNM-82,dz2}}. The function $E$  is called the
Cauchy kernel function.
}
\eex

\bex\label{ex 2.2}{
The hypercomplex variables
\beq \label{2.35}
z_j=z_j(w)=w_je_0-w_0e_j\,\,\,\, \big( j=1,\cdots,n\big)\eeq
are biregular \mbox{\rm \cite{BDS-RNM-82,Delanghe-70,Delanghe-72,Yeh-94,Yeh-91}}.
}
\eex
\bex\label{ex 2.3}
Let
$(\ell_1, \cdots,\ell_{k})\in N^{k}$,
i.e., $\ell_j$'s are $k$ elements out of $N$, where repetitions are allowed. We put
\beq\label{2.36}
V_{(\ell_1, \cdots,\ell_{k})}(w)
=\sum_{\pi(\ell_1,\cdots,\ell_k)}z_{\ell_1}(w)\,\cdots\,
z_{\ell_k}(w),\,\,\,\,\,w\in\mathbb{R} ^{n+1},
\eeq
where the sum runs over all the permutations of  $(\ell_1,\cdots,\ell_k)$, which is also the so-called the hypercomplex symmetric power. The hypercomplex symmetric power {$V_{(\ell_1, \cdots,\ell_{k})}$} is sometimes also
called the Fueter polynomial, being biregular\mbox{\rm\cite{BDS-RNM-82, TR}}.
\eex

\bex\label{ex 2.4}{
All derivatives  of $E$
\beq\label{2.37}
W_{\{\ell_1,\ell_2,\cdots\ell_k\}}(w)=(-1)^{k}
\displaystyle\frac{\partial^{k}\,E}{\partial w_{\ell_1}\partial w_{\ell_2}\cdots\partial w_{\ell_k}} (w),
\hspace{4mm}w\in\mathbb{ R}^{n+1}\setminus\{0\}
\eeq
are biregular\mbox{\rm \cite{BDS-RNM-82,Delanghe-70,Delanghe-72,Yeh-94,Yeh-91}}, being sometimes called negative powers \mbox{\rm \cite{TR}},  and
\beq\label{2.38}
W_{\{\ell_1,\ell_2,\cdots\ell_k\}}(w)=
O\Big(|w|^{-(n+k)}\Big)\,\,\,\mbox{ near}\,\,\,\,
\infty.
\eeq
  }
\eex

\mysection{Formulation of the {Hilbert BVPs} on $\mathbb{R}^{n+1}_{+}$}
To suitably present and formulize Hilbert boundary value problems on the Poincar\'{e} upper half space $\mathbb{R}^{n+1}_{+}$, we must introduce
a suitable
statement for the growth condition at
 the infinity for \vspace{0.5mm}the regular functions on the Poincar\'{e} upper half space $\mathbb{R}^{n+1}_{+}$.
 \vspace{2mm}

\subsection{Order at the infinity}

\mbox{}\hspace{5mm}
Assume $F$ is a regular
function on  the Poincar\'{e}  hyperplane $\mathbb{R}_{+}^{n+1}$ denoted as $\Phi\in
\mathcal {M}\big(\mathbb{R}^{n+1}_{+}\big)$.
We sometimes need  the concept of the order of $F$ at the infinity.
\bde\label{de 3.1}
Let $F\in
\mathcal {M}\left(\mathbb{R}^{n+1}_{+}\right)$, $m$ be a integer. If
\beq\label{3.1}
0<\beta=\limsup\limits_{w\in \mathbb{R}_{+}^{n+1},\, w\rightarrow\infty}\left|w^{-m}F(w)\right|<+\infty,
\eeq
then $F$ is said to be {of order $m$} at $w=\infty$, denoted it by $\mbox{\rm Ord\,}(F,\infty)=m$.
\ede

Sometimes it is more convenient to write
\beq
\label{3.2}
|F(w)|\stackrel{sup}{\approx}|G(w)|\,\,\,\mathrm{\it near}\,\,\infty\,
\,\,\,\,\,\mathrm{\it when}\,\,\,\,\,
  a\leq\limsup_{w\rightarrow\infty}\displaystyle\frac{|F(w)|}{|G(w)|}
 \leq A,
\eeq
where $F$ and $G$ $(|G(w)|>0)$ are two hypercomplex functions defined near $\infty$,  $A>a>0$ are two  constants. \vspace{2mm}

The following lemma is an obvious fact.
\ble\label{le 3.1}
$\big|F(w)\big|\stackrel{sup}{\approx}\big|w^{k}\big|$ near  $w=\infty$ if and only if\, $\mbox{\rm Ord\,}(F,\infty)=k$.
\ele

For  boundary behavior of function $\Phi\in
\mathcal {M}\big(\mathbb{R}^{n+1}_{+}\big)$ at the  infinity, there are {commonly three types of formulations:}
\beq\label{3.3}
\mbox{}\hspace{-18mm}(A)\,\,\, \lim\limits_{w\rightarrow\infty, w\in \mathbb{R} ^{n+1}_{+}}w^{-(m+1)}\Phi(w)=0, \,\,\,namely\,\,\,  \Phi(w)=o\left(w^{m+1}\right)\,\, near \,\,w=\infty,
\eeq
\beq\label{3.4}
\mbox{}\hspace{-18mm}(B)\,\,\,\limsup\limits_{w\in\mathbb{R}^{n+1}_{
+},\, w\rightarrow\infty}|w|^{-m}|\Phi(w)|=\beta,\,\,\, namely\,\,\, |\Phi(w)|=O\left(w^{m}\right)\,\,\, near\,\,\,w=\infty,
\eeq
\beq\label{3.5}
\mbox{}\hspace{-19mm}(C)\,\,\,\,\,\mbox{\rm Ord\,}(\Phi,\infty)\leq m, \,\,\,namely \,\,\, \,\,\, |\Phi(w)|\stackrel{sup}{\approx}\big|w^{k}\big|\,\,\,\mathrm{\it near }\,\,\, w=\infty\,\,\,\,{\it when} \,\,\,k\leq m.
\eeq

\bre\label{re 3.1}
{\rm
Obviously,  (C) implies (B), while (B) implies (A). So, the condition (A) is the weakest.
In \cite{ddq17}, We {used} the condition (A) in  Riemann  boundary value problems with the hyperplane as  jump surface, which {was} an innovation.
In the present paper, {we will still use it in our} Hilbert  boundary value problems on the   Poincar\'{e} upper  half space.
}
\ere

Some symbols will be  used in the following \cite{Yeh-94}.
Let
\beq\label{3.6}
Z(w)=\big(z_1(w),z_2(w),\cdots,z_n(w)\big)\,\,\,\,\mathrm{\it and }\,\,\,\,\alpha=\big[\alpha_1,\alpha_2,\cdots,\alpha_n\big],
\eeq
 where
$z_{j}$'s are
the hypercomplex variables given in Example \ref{ex 2.2}
and $\alpha_j$'s are nonnegative integers. Then the symmetry power $Z^\alpha$ in $\mathbb{R}^{n+1}$ is defined as the sum of all {possible} $z_i$ products of which each {contains $z_i$ factor exactly $\alpha_i$ times}.
For example, for $n=2$
\beq\label{3.7}
(z_1,z_2)^{[0,0]}=1,\,\,\,
(z_1,z_2)^{[1,1]}=z_1z_2+z_2z_1,\,\,\,(z_1,z_2)^{[2,0]}=z_1^2.
\eeq

Introduce the mapping
\beq\label{3.8}
\mbox{\bf \LARGE $ \alpha$} :\big(\ell_1,\ell_2,\cdots,\ell_k\big)\mapsto \alpha=\big[\alpha_1,\cdots\alpha_n\big],
\eeq
{where  $\alpha_j$ is  the number of times that $j$ appears in $\big(\ell_1,\ell_2,\cdots,\ell_k\big)\in N^{k}$ and  $N=\{1,...,n\}$}.
Now we may  rewrite $W_{\{\ell_1,\ell_2,\cdots,\ell_k\}}$ in  Example \ref{ex 2.4} by
\beq \label{3.9}
(-1)^{|\alpha|} W_{\{\ell_1,\ell_2,\cdots,\ell_k\}}(w)
=\displaystyle
\frac{\partial^{|\alpha|}\,E}{\partial w_0^{\alpha_0}\cdots\partial w_n^{\alpha_n}} (w)
\equiv\displaystyle
\frac{\partial^{|\alpha|}\,E}{\partial w^{\alpha}} (w)
\equiv
\Big(\partial^{\alpha}[E]\Big)(w)
,
\,\,w\in\mathbb{R}^{n+1}\backslash\{0\}.
\eeq
So, $w^{n+|\alpha|}\big(\partial^{\alpha}[E]\big)(w)$ is bounded {according to} Example \ref{ex 2.4}.

Similarly,
\beq\label{3.10}
Z^{\alpha}(w)=V_{(\ell_1, \cdots,\ell_{k})}(w)
=\sum_{\pi(\ell_1,\cdots,\ell_k)}z_{\ell_1}\,\cdots\,z_{\ell_k}(w)\,\,\,\,\mathrm{\it with}\,\,\,\, k=|\alpha|,
\eeq
where $V_{(\ell_1, \cdots,\ell_{k})}$ is the Fueter polynomial given in Example \ref{ex 2.3}.
So $Z^{\alpha}$ is a
biregular function
\cite{BDS-RNM-82}.

\bde\label{de 3.2}
Let $m\geq0$. We call
\beq\label{3.11}
f(w)=\sum_{|\alpha|=0}^{m} Z^{\alpha}(w)\,c_{\alpha}\,\, \,\,\mbox{with}\,\,\,\,
\sum_{|\alpha|=m}|c_{\alpha}|\not=0
\eeq
 a hypercomplex symmetric  polynomial of degree $m$. In such case we denote $\mbox{\rm Deg}(f)=m$.
\ede

\ble [\mbox{\rm see \cite{ld}}]\label{le 3.2}
Let $f$ be a hypercomplex symmetric polynomial. Then  $\mbox{\rm Deg}(f)=m$ if and only if $\mbox{\rm Ord}(f,\infty)=m$.
\ele

\subsection{$\widehat{H}$ class of functions}

\vspace{1mm}\mbox{}\hspace{6mm}In order to state the condition of the input function for {Hilbert BVPs}, we need to introduce some classes of hypercomplex functions used frequently in this paper.

 \bde  \label{de 3.3}
  Assume $f$ is defined on $\Omega\subseteq\mathbb{R}^{n+1}$.
If\vspace{-1mm}
\beq \label{3.12}
\big|f(t)-f(s)\big|\leq M\big|t-s\big|^{\mu}\,\,\,\,\,\,(0<\mu\leq1)
 \eeq
 for arbitrary points
$t,\,s$ on $\Omega$, where $M$ and $\mu$ are  constants, then
$f$ is said to satisfy  H\"{o}lder condition of order $\mu$,
denoted by $f\in H^{\mu}(\Omega)$. The constants  $\mu$ and $M$ are called, respectively,  the H\"{o}lder
index and {the}  H\"{o}lder coefficient of $f$. If the order $\mu$ is not emphasized, it may be denoted
briefly {as} $f\in  H(\Omega)$.
\ede

  Sometimes we will use $\displaystyle\frac{1}{w}$ \vspace{1.5mm}to represent $w^{-1}$ for $w\in \mathbb{R}^{n+1}$,
which suggests the similarity with some results in the classical complex analysis.

 \bde \label{de 3.4}
Assume $f$ is defined on $\Omega\subseteq\mathbb{R}^{n+1}$.
If
\beq \label{3.13}
\Big|f(\xi)-f(\zeta)\Big|\leq
M\left|\displaystyle\frac{1}{\xi}-\frac{1}{\zeta}\right|^{\mu} 
\,\,\,\,(0<\mu\leq1)
 \eeq
 for arbitrary points
$\xi,\, \zeta$ on $\Omega\setminus\{0\}$, where $M$ and $\mu$ are constants,
then
$f$ is said to satisfy  $H_{\dag}$ condition of order $\mu$,
 denoted by $f\in {H}_{\dag}^{\mu}\left(\Omega\right)$. The constants  $\mu$ and $M$ are called, respectively,  the $\dag$-H\"{o}lder
index and {the} $\dag$-H\"{o}lder coefficient of $f$. If the order $\mu$ is not emphasized, it may be denoted
briefly {as} $f\in H_{\dag}(\Omega)$.
\ede

\bde \label{de 3.5}
If $f\in H^{\mu}(\Omega)\cap H_{\dag}^{\mu}(\Omega)$, \vspace{1.5mm} then
$f$ is said to satisfy the  $\widehat{H}$ condition of order $\mu$ on $\Omega$,
 denoted by $f\in \widehat{H}^{\mu}\left(\Omega\right)$ or briefly $f\in \widehat{H}(\Omega)$.
\ede

The conditions $(\ref{3.12})$ and $(\ref{3.13})$ are, respectively, called the H\"{o}lder condition and $\dag$-H\"{o}lder condition of the $\widehat{H}(\Omega)$ class function $f$.
More discussions of the two classes of functions   may be found in \cite{ddq17}.
\bde\label{de 3.6}
Let $f$ be a function defined on   $\Omega$ with $\infty$ as its  cluster point.
If \vspace{-2mm}
\beq
\label{3.14}
f(\infty)=\lim\limits_{w\in\Omega,\,w\rightarrow\infty}f(w)\eeq
  exists and
\beq\label{3.15}
\Big|f(w)-f(\infty)\Big|\leq \frac{M}{|w|^{\mu}}\,\,\Big(0<\mu\leq1\Big),\,\,\,\,\,w\in\Omega\backslash\{0\},
\eeq
where $M$ is a constant, then we say that $f$ satisfies the pointwise $\dag$-H\"{o}lder condition at the infinity in $\Omega$,
denoted by $f\in H_{\dag}^{\mu}(\Omega,\infty)$, or also briefly  $f\in H_{\dag}^{\mu}(\infty)$ and $f\in H_{\dag}(\infty)$.
\ede

In the sequel the following notations will be used.
Let
\beq \label{3.16}
f_{\mathbbm{m}}(w)=w^{m} f(w).
\eeq

$(1)$ If  $f_{\mathbbm{m}}\!\in\! {H}^{\mu}\left(\mathbb{R} ^{n+1}_{0}\right)$,
  then we write  $f\!\in\! {H}_{m}^{\mu}\big(\mathbb{R}^{n+1}_{0}\big)$, or briefly  $f\in H_{m}\big(\mathbb{R} ^{n+1}_{0}\big)$.\vspace{1mm}\vspace{1mm}

  $(2)$ If  $f_{\mathbbm{m}}\in {H}_{\dag}^{\mu}\big(\mathbb{R}^{n+1}_{0}\big)$,\vspace{1mm}
  then we write   $f\in {H}_{m, \dag}^{\mu}\big(\mathbb{R}^{n+1}_{0}\big)$,
  or briefly $f\in H_{m,\dag}\big(\mathbb{R}^{n+1}_{0}\big)$.\vspace{1mm}

   $(3)$ If  $f_{\mathbbm{m}}\in {H}_{\dag}^{\mu}\big(\mathbb{R}^{n+1}_{0},\infty\big)$,\vspace{1mm}
  then we write   $f\in {H}_{m, \dag}^{\mu}\big(\mathbb{R}^{n+1}_{0},\infty\big)$, or briefly   $f\in H_{m,\dag}\big(\infty\big)$.

The following classes of functions will also be used in Hilbert boundary value problems, of which the details can be found in \cite{ddq17}:
\beq\label{3.17}
\widehat{H}_{m}\Big(\mathbb{R} ^{n+1}_{0}\Big)= {H}_{m}\Big(\mathbb{R} ^{n+1}_{0}\Big)\bigcap {H}_{m, \dag}\Big(\mathbb{R} ^{n+1}_{0}\Big),
\eeq
\beq\label{3.18}
\widehat{H}_{m,0}\Big(\mathbb{R} ^{n+1}_{0}\Big)= \widehat{H}_{m}\Big(\mathbb{R} ^{n+1}_{0}\Big)\bigcap\Big\{f, f_{\mathbbm{m}}(\infty)=0\Big\}.
\eeq

\subsection{Formulation of {Hilbert boundary value problems}}
\mbox{}\hspace{6mm}
Let $\Phi$ be a regular function defined on the Poincar\'{e} upper half space
$\mathbb{R} _{+}^{n+1}$ which can extend continuously to the hyperplane
$\mathbb{R} _{0}^{n+1}$. For clarity, we denote  its boundary value by

\beq\label{3.19}
\Phi^{+}(t)=\lim\limits_{w\rightarrow t,\,
w\in\mathbb{R} _{+}^{n+1}}\Phi(w),\,\,\,\,\,t\in\mathbb{R} _{0}^{n+1}.
\eeq

\mbox{\bf The Hilbert boundary value problem (or simply $\mathrm{H}_{m}$ problem).}
 \,\,Find a
left (right) regular function in $\mathbb{R} _{+}^{n+1}$ which can extend
continuously to $\mathbb{R} _{0}^{n+1}$ such that
\beq \label{3.20}\left\{
\begin{array}{l}
\big(D[\Phi]\big)(w)=0\,\,\,\,\mbox{\it for}\,\,\,\,w\in \mathbb{R}^{n+1}_{+}\,\,\,(\mbox{\it regularity}),\\[2mm]
\mbox{\rm Re}\Big\{\Phi^{+}(x)\,\lambda\Big\} = c(x)\,\,\,\mbox{\it for}\,\,\, x\in\mathbb{R}_{0}^{n+1}\,\,\,
(\mbox{\it boundary condition}),\\[2mm]
\Phi(w)=o\big(w^{m+1}\big)\,\,\,\mbox{\it near}\,\,\,\infty\,\,\,
(\mbox{\it growth condition} ),
\end{array}\right.
\eeq
where $\lambda$ is a given  constant whose inverse  exists, $c$ is a given
para real-valued function and

\beq \label{3.21}
c\in\left\{
\begin{array}{ll}
\widehat{H}\Big(\mathbb{R}^{n+1}_{0}\Big),&\mathrm{\it when}\,\,\,\,\,
m\geq0\,\,\,\mathit{(non-negative order)},\\[2mm]
\widehat{H}_{r,0}\Big(\mathbb{R}^{n+1}_{0}\Big),&\mathrm{\it when}\,\,\,\,\,
m<0\,\,\,\,\,\mbox{\it with}\,\,\,\,\,r=-(m+1).
\end{array}
\right.
\eeq

Obviously, the above BVPs are the direct generalizations of the classical
 Hilbert BVPs \cite{Lu,Mu} from complex plane to the space $C(V_{n})$.
In the classical complex analysis, $\lambda$ is allowed to be a complex function                              when
$\Phi$ is finite at the infinity, $i,e$., when $m=0$.

A great number of literatures have discussed the solution methods of  BVPs for some special case  in
Clifford analysis. One methodology is the so-called symmetric extension of which  the main idea is to translate Hilbert BVPs (\ref{3.20}) to  {the equivalent} Riemann problems with {some} additional restriction condition. In the complex analysis, this approach has been  quite successful \cite{Lu, Mu}, but in the higher dimensional spaces {there} exist many
obstacles to implement  this approach. Not all BVPs can be solved by the symmetric extension, especially
for the case when $\lambda$  in (\ref{3.20}) is a hypercomplex function. {For} the case when $\lambda$ is a hypercomplex function, it can still be  translated to a Riemann problem, but an open problem.
 This method
{was} tried by Xu and Zhou in \cite{Xu} with $m=0$.  {Gong and Du solved the $\mathrm{H}_{m}$ problems with $m\geq0$ in \cite{GD} by using the symmetric extension method that generalizes the results in \cite{Lu, Mu} directly to the Clifford setting}.

The simplest {$\mathrm{H}_{m}$ problem} is the Schwarz problem, $i.e.$, $\lambda=1$.
Let us start with some {discussions} on the Schwarz problem. \vspace{1mm}

\mbox{\bf The Schwatz boundary value problem $\mathrm{S}_{m}$.}
 \,\,Find a
left (right)  regular function $\Phi$ in $\mathbb{R} _{+}^{n+1}$ {that} can extend
continuously to $\mathbb{R} _{0}^{n+1}${such that}
\beq \label{3.22}\left\{
\begin{array}{l}
\big(D[\Phi]\big)(w)=0\,\,\,\,\mbox{\it  for}\,\,\,\,w\in \mathbb{R}^{n+1}_{+}\,\,\,\,(\mbox{\it regularity}),\\[2mm]
\mbox{\rm Re}\Big\{\Phi^{+}(x)\Big\} = c(x)\,\,\,\mbox{\it for}\,\,\, x\in\mathbb{R}_{0}^{n+1}\,\,\,
(\mbox{\it boundary condition}),\\[2mm]
\Phi(w)=o\big(w^{m+1}\big)\,\,\,\mbox{\it near}\,\,\,\infty\,\,\,
(\mbox{\it growth condition} ),
\end{array}\right.
\eeq
where $c$ is a given
para real-valued function and $c$ satisfies $(\ref{3.21})$. \bre{\rm Below, we generally assume that $\Phi$ is left regular.}
\ere

\mysection{Symmetric extension and self-reflex action}

In order to solve the $\mathrm{H}_m$ problem, we try to transfer it
into an Jump problem {$R_{m}$}  discussed in  \S 6.3 of \cite{ddq17}.
To do so, we introduce the symmetric extension and the self-reflex action.
\vspace{1.5mm}

\subsection{Symmetric extension}
\mbox{}\hspace{6mm}
Introduce the reflection operator, $*:C(V_n)\longrightarrow C(V_n)$, with respect to the hyperplane $\mathbb{R}_{0}^{n+1}$, by
\beq \label{4.1}
\lambda^{*}=\mbox{\rm Re} (\lambda)-e_{n}\,\mbox{\rm Im}^l (\lambda)\,,\,\,\,\lambda\in C(V_n),
\eeq
or
\beq\label{4.2}
\lambda^{*}=x-e_n\,y^{l},\,\,\,x,y^l\in C\big(V_{n-1}\big),
\eeq
where
\beq\label{4.3}
x={\mbox{\rm Re} (\lambda)},\,\,\,\,y^l=\mbox{\rm Im}^l (\lambda).
\eeq
Obviously, by  (\ref{2.15})
\beq\label{4.4}
\lambda^*=x-y^{r}\,e_n,\,\,\,\,\,\lambda\in C(V_n),
\eeq
where
\beq\label{4.5}
x=\mbox{\rm Re} (\lambda),\,\,\,\,y^r=\mbox{\rm Im}^r (\lambda).
\eeq
In particular,
\beq\label{4.6}
w^{*}=\sum_{j=0}^{n-1}w_j e_j-w_n e_n\,\,\,\,\,\,\mathrm{\it when}\,\,\,\,\,\,
w=\sum_{j=0}^{n}w_j e_j
=\mbox{\rm Re} (w)+\mbox{\rm Im} (w)\,e_n\in \mathbb{R}^{n+1}.
 \eeq
 We know
that $\lambda^{*}$ and $\lambda$ are a pair of points symmetric to the hyperplane $\mathbb{R}_{0}^{n+1}$ and
\beq\label{4.7}
\big(\lambda^{*}\big)^{*}=\lambda,\,\,\,\,\lambda\in C(V_n).
\eeq
Moreover, by (\ref{2.3}) and (\ref{4.1}),
\beq\label{4.8}
\big(\overline{w}\big)^{*}=\overline{w^{*}}, \,\,\,\,\,w\in\mathbb{R}^{n+1}.
\eeq

For a function $\Phi(w)$ in $\mathbb{R}_{+}^{n+1}$, we define a function
 on $\mathbb{R}^{n+1}_{-}$ by
\beq
\label{4.9}
{\Phi}^{*}(w)=\big[\Phi\left(w^{*}\right)\big]^{*},
\,\,\,\,\,\,\,w\in\mathbb{R}_{-}^{n+1},
\eeq
which is called its accompanying function. More specifically,
if
\beq\label{4.10}
{\Phi}\big(w\big)=\big(\mathrm{Re}\,\Phi\big)\big(w\big)
+e_n\,\big(\mathrm{Im}^l\,\Phi\big)\big(w\big),\,\,\,\,\,w\in\mathbb{R}^{n+1}_{+},
\eeq
with
\beq\label{4.11}
\big(\mathrm{Re}\,\Phi\big)\big(w\big)=\displaystyle\sum_{A\in\mathcal {P} \{1,\cdots,n-1\}}u_{A}\big(w_0,\cdots,w_{n-1},w_n\big)e_{A}
\eeq
and
\beq\label{4.12}
\big(\mathrm{Im}^l\,\Phi\big)\big(w\big)=\displaystyle\sum_{A\in\mathcal {P} \{1,\cdots,n-1\}}v_{A}\big(w_0,\cdots,w_{n-1},w_n\big)e_{A},
\eeq
where $u_A$ and $v_A$ are para real-valued functions,
then
\beq\label{4.13}
{\Phi}^{*}\big(w\big)=
\big(\mathrm{Re}\,\Phi^{*}\big)\big(w\big)
+e_n\,\big(\mathrm{Im}^{l}\,\Phi^{*}\big)\big(w\big),
\,\,\,\,\,\,w\in\mathbb{R}^{n+1}_{-},
\eeq
where
\beq\label{4.14}
\big(\mathrm{Re}\,\Phi^{*}\big)\big(w\big)=
\displaystyle\sum_{A\in\mathcal {P} \{1,\cdots,n-1\}}u_{A}\big(w^{*}\big)e_{A}=\displaystyle\sum_{A\in\mathcal {P} \{1,\cdots,n-1\}}u_{A}\big(w_0,\cdots,w_{n-1},-w_n\big)e_{A},
\eeq
\beq\label{4.15}
\big(\mathrm{Im}^{l}\,\Phi^{*}\big)\big(w\big)
=-
\displaystyle\sum_{A\in\mathcal {P} \{1,\cdots,n-1\}}v_{A}\big(w^{*}\big)e_{A}
=-\displaystyle\sum_{A\in\mathcal {P} \{1,\cdots,n-1\}}v_{A}\big(w_0,\cdots,w_{n-1},-w_n
\big)e_{A}.
\eeq

Similarly, if the original function $\Phi(w)$ is defined in $\mathbb{R}_{-}^{n+1}$, then its accompanying function
  ${\Phi}^{*}(w)$ is determined by
 \beq
\label{4.16}
{\Phi}^{*}(w)=\big[\Phi\left(w^{*}\right)\big]^{*},\,\,\,\,\,\,\,w\in\mathbb{R}_{+}^{n+1}.
\eeq

We easily see that
\beq \label{4.17}
\big(\Phi^{*}\big)^{*}(w)=\Phi(w),\,\,\,\,\,w\in\mathbb{R}_{+}^{n+1}\,\,
\Big(\mathbb{R}_{-}^{n+1}\Big),
\eeq
that is to say, if $\Phi(w)$ is acted by the reflection oprator  with respect to $\mathbb{R}_{0}^{n+1}$ twice, then it returns to $\Phi(w)$ itself, or the reflection operator $*$ is idempotent.

Obviously, when
(\ref{4.11}) and (\ref{4.12}) hold for $w\in\mathbb{R}_{-}^{n+1}$ then (\ref{4.14}) and (\ref{4.15}) also hold for $w\in\mathbb{R}_{+}^{n+1}$.
\bre \label{re 4.1}
{\rm When the original function $\Phi$ is defined on $\mathbb{R}^{n+1}_{+}$ or
$\mathbb{R}^{n+1}_{-}$, we have that
\beq\label{4.18}
\begin{array}{ll}
&\Phi(w)\hspace{4mm} \Big(w=\big(w_0,w_1,\cdots,w_{n}\big)\in \mathbb{R}^{n+1}_{+}\,\,
\Big(
\mathbb{R}^{n+1}_{-}\Big)\Big)\\[4mm]
=&\!\!\!\!\displaystyle\sum_{A\in\mathcal {P} \{1,\cdots,n-1\}}u_{A}\big(w_0,\cdots,w_{n-1},w_n\big)e_{A}
\,+e_n\,\displaystyle\sum_{A\in\mathcal {P} \{1,\cdots,n-1\}}v_{A}\big(w_0,\cdots,w_{n-1},w_n\big)e_{A}
\end{array}
\eeq
is equivalent to
\beq\label{4.19}
\begin{array}{ll}
&\Phi^{*}(w)\hspace{4mm} \Big(w=\big(w_0,w_1,\cdots,w_{n}\big)\in \mathbb{R}^{n+1}_{-}\,\,\Big(
\mathbb{R}^{n+1}_{+}\Big)\Big)\\[4mm]
=&\displaystyle\!\sum_{A\in\mathcal {P} \{1,\cdots,n-1\}}\!\!u_{A}\big(w^{*}\big)e_{A}
-e_n\,\displaystyle\!\sum_{A\in\mathcal {P} \{1,\cdots,n-1\}}v_{A}\!\big(w^{*}\big)e_{A}\\[9mm]
=&\displaystyle\!\sum_{A\in\mathcal {P} \{1,\cdots,n-1\}}\!u_{A}\big(w_0,\cdots,w_{n-1},-w_n\big)e_{A}
\!-\!\displaystyle\!\sum_{A\in\mathcal {P} \{1,\cdots,n-1\}}\!v_{A}\big(w_0,\cdots,w_{n-1},-w_n\!\big)e_{A}e_n.
\end{array}
\eeq
 }
 \ere

 \ble[\!\mbox{\rm\cite{ddq17}}]\label{le 4.1}
Let $f(w)$ and $g(w)$ be hypercomplex functions defined on $\mathbb{R} _{+}^{n+1}$ $\Big(\mathbb{R} _{-}^{n+1}\Big)$, then
\beq\label{4.20}
\Big|f^{*}\big(w^*\big)\Big|=
\Big|f(w)\Big|,\,\,\,\,x\in\mathbb{R} _{+}^{n+1}\,\,\,\,\,\,\Big(\mathbb{R} _{-}^{n+1}\Big),
\eeq
as well as
\beq\label{4.21}
\Big[f(w)g(w)\Big]^{*}=f^{*}\big(w^{*}\big)g^{*}\big(w^{*}\big),\,\,\,\,w\in\mathbb{R} _{+}^{n+1}\,\,\,\,\,\Big(\mathbb{R} _{-}^{n+1}\Big).
\eeq

If the inverse function $f^{-1}$ exists, then
the inverse function of $f^{*}$ also exists,  and
\beq\label{4.22}
\Big[f^{*}\Big]^{-1}
=\Big[f^{-1}\Big]^{*},
\eeq
more precisely,
\beq\label{4.23}
\Big[f^{*}\Big]^{-1}\big(s^{*}\big)
=\Big[f^{-1}\Big]^{*}\big(s^{*}\big),\,\,\,\,s
\in f\Big(\mathbb{R} _{+}^{n+1}\Big)\,\,\,\,\bigg(f\Big(\mathbb{R} _{-}^{n+1}\Big)\bigg),
\eeq

\ele

\p $(\ref{4.20})$ comes directly from $(\ref{4.18})$ and $(\ref{4.19})$ by $(\ref{2.16})$.   Just to be clear, in (\ref{2.20}) we rewrite $u_A$ and $v_A$ {as, respectively,} $u_{f,A}$ and $v^l_{f,A}$.
Then,
\beq\label{4.24}
\begin{array}{ll}
&\Big[f(w)\,g(w)\Big]^{*}\,\,\,\,\,\Big( \mathrm{\it say}\,\,\,w\in\mathbb{R} _{+}^{n+1}\Big)\\[6mm]
=&
\displaystyle\left[\sum_{A,B\in\mathcal {P} \{1,\cdots,n-1\}}
\Big\{u_{f,A}(w)e_{A}+e_n\,v^{l}_{f,A}(w)e_{A}\Big\}\Big\{
u_{g,B}e_{B}+e_n\,v^{l}_{g,B}(w)
e_{B}\Big\}\right]^{*}\\[8mm]
=&\displaystyle\sum_{A,B\in\mathcal {P} \{1,\cdots,n-1\}}
\left[u_{f,A}(w)u_{g,B}(w)
-(-1)^{\#(A)}v^{l}_{f,A}(w)v^{l}_{g,B}(w)
\right]e_{A}e_{B}\\[8mm]
&\,\,-\,e_n\,
\displaystyle\sum_{A,B\in\mathcal {P} \{1,\cdots,n-1\}}
\!\left[v^{l}_{f,A}(w)u_{g,B}(w)
-(-1)^{\#(A)}u_{f,A}(w)v^{l}_{g,B}(w)
\right]e_{A}e_{B}.
\end{array}
\eeq
On the other hand, by (\ref{4.11}), (\ref{4.12}), (\ref{4.14}),
(\ref{4.15}) and (\ref{4.2}), we have
\beq\label{4.25}
\mathrm{Re}\Big[f^{*}g^{*}\Big]\big(w^{*}\big)=
\displaystyle\sum_{A,B\in\mathcal {P} \{1,\cdots,n-1\}}
\left[u_{f,A}(w)u_{g,B}(w)
-(-1)^{\#(A)}v^{l}_{f,A}(w)v^{l}_{g,B}(w)
\right]e_{A}e_{B},
\eeq
and
\beq\label{4.26}
\mathrm{Im}^{l}\Big[f^{*}g^{*}\Big]\big(w^{*}\big)=
-\displaystyle\sum_{A,B\in\mathcal {P} \{1,\cdots,n-1\}}
\!\left[(-1)^{\#(A)}u_{f,A}(w)v^{l}_{g,B}(w)
+v^{l}_{f,A}(x)u_{g,B}(w)
\right]e_{A}e_{B}.
\eeq

(\ref{4.24}), (\ref{4.25})  and (\ref{4.26}) imply (\ref{4.21}).

 Noting that
 \beq\label{4.27}
 w=\Big(f^{-1} \circ f\Big)(w)=f^{-1}\big(f(w)\big),\,\,\,\,\,w\in\mathbb{R} _{+}^{n+1}\,\,\Big(\mathbb{R} _{-}^{n+1}\Big),
 \eeq
we have, by $(\ref{4.7})$  and $(\ref{4.9})$,
\beq\label{4.28}
\Big(\big[f^{-1}\big]^{*} \circ f^{*}\Big)\big(w^{*}\big)
=\Big[f^{-1}\Big(\big[f^{*}\big]^{*}\big( \left[w^{*}\right]^{*}\big)\Big)\Big]^{*}
=\Big[f^{-1}\big(f(w)\big)\Big]^{*}
=w^{*},\,\,\,w^{*}\in\mathbb{R} _{-}^{n+1}\,\,\Big(\mathbb{R} _{+}^{n+1}\Big),
\eeq
which is (\ref{4.22}).
\ef

Introduce the operators
\beq\label{4.29}
\displaystyle\frac{\partial}{\partial x}=
\sum_{j=0}^{n-1}e_j\,\frac{\partial}{\partial w_j}\,\,\,\,\,\,\mathrm{\it with}\,\,\,\,\,\,x=\mathrm{Re}(w)
\eeq
and
\beq\label{4.30}
\displaystyle\frac{\partial}{\partial y}=
\frac{\partial}{\partial w_n}\,\,\,\,\,\mathrm{\it with}\,\,\,\,\,y=\mathrm{Im}(w).
\eeq
Their actions on functions from the left and from the right {are} governed by the rules
\beq\label{4.31}
\frac{\partial [f]}{\partial x}=
\sum_{k=0}^{n-1}\sum_{A}e_k\,e_{A}\frac{\partial f_A}{\partial w_k},\,\,\,\,\,\,
\frac{[f] \partial}{\partial x}=
\sum_{k=0}^{n-1}\sum_{A}e_A\,e_{k}\frac{\partial f_A}{\partial w_k}.
\eeq

\ble[{\bf Cauchy-Riemann equations} {\rm \cite{Xu, GD}}]\label{le 4.2}

$(1)$\, $\Phi$ is left regular on $\mathbb{R}_{+}^{n+1}$ $\Big(\mathbb{R}_{-}^{n+1}\!\Big)$ if and only if
\beq\label{4.32}
\displaystyle\frac{\partial [U]}{\partial x}(w)=\displaystyle\frac{\partial \big[V^{l}\big] }{\partial y}(w)
\,\,\,\,\,\mathrm{\it and}\,\,\,\,\,\,
\displaystyle\frac{\partial [U]}{\partial y}(w)=-\displaystyle\frac{\partial \big[V^{l}\big] }{\partial x}(w),
\eeq
where
\beq\label{4.33}
U(w)=\big(\mathrm{Re}\,\Phi\big)\big(w\big)\,\,\,\,\,\mathrm{\it and}\,\,\,\,\,\,
V^{l}(w)=\big(\mathrm{Im}^{l}\,\Phi\big)\big(w\big)
\eeq
{are, respectively,} the real and the left imaginary parts of $\Phi$. \vspace{1mm}

$(2)$\,  $\Phi$ is right regular on $\mathbb{R}_{+}^{n+1}$ $\Big(\mathbb{R}_{-}^{n+1}\Big)$ if and only if
\beq\label{4.34}
\displaystyle\frac{[U]\,\partial }{\partial x}(w)=\displaystyle\frac{\big[V^{r}\big] \partial}{\partial y}(w)
\,\,\,\,\,\mathrm{\it and}\,\,\,\,\,
\displaystyle\frac{[U]\,\partial }{\partial y}(w)=-\displaystyle\frac{\big[V^{r}\big]\partial }{\partial x}(w),
\eeq
where
\beq\label{4.35}
U(w)=\big(\mathrm{Re}\,\Phi\big)\big(w\big)\,\,\,\,\,\mathrm{\it and}\,\,\,\,\,
V^{r}(w)=\big(\mathrm{Im}^{r}\,\Phi\big)\big(w\big)
\eeq
{are, respectively,} the real and the right imaginary parts of $\Phi$.
\ele

\p\,If $\Phi$ is left regular on $\mathbb{R}_{+}^{n+1}$, it is easy to see, by Remark 4.1, that
\beq\label{4.36}
\begin{array}{lll}
&&\Big(D[\Phi]\Big)(w) \,\,\,\Big(w\in\mathbb{R}_{+}^{n+1}\Big)\\[6mm]
&=&\displaystyle \sum_{k=0}^{n-1}e_k\left[\sum_{A\in\mathcal {P} \{1,\cdots,n-1\}}\, e_{k}\,e_A\frac{\partial u_A}{\partial w_k}\big(w\big)
+e_n\,\sum_{A\in\mathcal {P} \{1,\cdots,n-1\}}e_A\frac{\partial v^{l}_A}{\partial w_k}\big(w\big)\right]
\\[9mm]
 &&+\displaystyle\,
 e_n\left[\sum_{A\in\mathcal {P} \{1,\cdots,n-1\}}e_A \frac{\partial u_A}{\partial w_n}\big(w\big)
 +e_n\,\sum_{A\in\mathcal {P} \{1,\cdots,n-1\}}e_{A}\frac{\partial v^{l}_A}{\partial w_n}\big(w\big)\right]\,\,\,\,\,\,\,
 \Big(\mathrm{\it by}\,\,\,\,
(\ref{4.10})\Big)\\[8mm]
&=&\displaystyle \sum_{k=0}^{n-1} \sum_{A\in\mathcal {P} \{1,\cdots,n-1\}}\, e_{k}\,e_A\frac{\partial u_A}{\partial w_k}\big(w\big)
-\sum_{k=0}^{n-1} \sum_{A\in\mathcal {P} \{1,\cdots,n-1\}}e_A\frac{\partial v^{l}_A}{\partial w_k}\big(w\big)
\\[8mm]
 &&+\displaystyle\sum_{A\in\mathcal {P} \{1,\cdots,n-1\}}e_n\,e_A \frac{\partial u_A}{\partial w_n}\big(w\big)
 +\sum_{A\in\mathcal {P} \{1,\cdots,n-1\}}e_n\,e_A \frac{\partial v^{l}_A}{\partial w_n}\big(w\big)\,\,\,\,\,\,\,
 \\[8mm]
 &=&\displaystyle
 \frac{\partial U}{\partial x}(w)-
 \displaystyle\frac{ \partial\big[V^{l}\big]}{\partial y}(w)
+e_n\left[
\displaystyle\frac{\partial U}{\partial y}(w)+\displaystyle\frac{\partial \big[V^{l}\big]}{\partial x}(w)\right]\,\,\,\Big(\mathit{ by }\,\,\,(\ref{4.29})
\,\,\,and\,\,\,(\ref{4.30}) \Big),
 \end{array}
 \eeq
 which implies $(1)$. The proof of $(2)$ is similar.
 \ef

\bth[Symmetry principle for regular functions]\label{th 4.1}
If $\Phi$ is regular on $\mathbb{R}_{+}^{n+1}\,\Big(\mathbb{R}_{-}^{n+1}\Big)$ then
$\Phi^{*}$ is also regular on $\mathbb{R}_{-}^{n+1}\,\Big(\mathbb{R}_{+}^{n+1}\Big)$.
\eeth

\p By (\ref{4.18}) and (\ref{4.19}), we have
\beq\label{4.37}
U^{*}(w)=\mathrm{Re}\,\Phi^{*}(w)=\mathrm{Re}\,\Phi(w^*)=U\big(w^*\big),\,\,\,\,
w\in\mathbb{R}_{-}^{n+1}\,\Big(\mathbb{R}_{+}^{n+1}\Big),
\eeq
and
\beq\label{4.38}
\big(V^{*}\big)^{l}\big(w\big)=\mathrm{Im}\,\Phi^{*}(w)
=-\mathrm{Im}\,\Phi\big(w^*\big)=-V^{l}\big(w^*\big),\,\,\,\,
w\in\mathbb{R}_{-}^{n+1}\,\Big(\mathbb{R}_{+}^{n+1}\Big).
\eeq
 So,
 \beq\label{4.39}
 \displaystyle\frac{\partial \big[U^*\big]}{\partial x}\big(w\big)=\displaystyle \sum_{k=0}^{n-1}e_k\sum_{A\in\mathcal {P} \{1,\cdots,n-1\}}e_A\frac{\partial u_A}{\partial w_k}\big(w_0,\cdots,-w_n\big)
 =\displaystyle\frac{\partial [U]}{\partial x}\big(w_0,\cdots,-w_n\big),
 \eeq
\beq\label{4.40}
 \displaystyle\frac{\partial \big[U^*\big]}{\partial y}\big(w\big)=\displaystyle \sum_{A\in\mathcal {P} \{1,\cdots,n-1\}}e_A\frac{\partial u_A}{\partial w_n}\big(w_0,\cdots,-w_n\big)
 =\displaystyle\frac{\partial [U]}{\partial y}\big(w_0,\cdots,-w_n\big),
 \eeq
 and
\beq\label{4.41}
\displaystyle\frac{\partial \big[(V^*)^{l}\big]}{\partial y}\big(w\big)
=\displaystyle \sum_{k=0}^{n-1}e_k\sum_{A\in\mathcal {P} \{1,\cdots,n-1\}}e_A\frac{\partial v^{l}_A}{\partial w_k}\big(w_0,\cdots,-w_n\big)
=
\displaystyle\displaystyle\frac{\partial \big[V^l\big]}{\partial y}\big(w\big),
\eeq
\beq\label{4.42}
\displaystyle\frac{\partial \big[(V^*)^{l}\big]}{\partial x}\big(w\big)
=\displaystyle \sum_{k=0}^{n-1}e_k\sum_{A\in\mathcal {P} \{1,\cdots,n-1\}}e_A\frac{\partial v^{l}_A}{\partial w_k}\big(w_0,\cdots,-w_n\big)
=
\displaystyle\displaystyle\frac{\partial \big[V^l\big]}{\partial x}\big(w\big).
\eeq

By Theorem \ref{th 4.1}, (\ref{4.39}), (\ref{4.40}), (\ref{4.41}), (\ref{4.42}),
we know that
\beq\label{4.43}
\Phi \,\,\mathrm{\it is\,\, (left)\,\,regular\,\,on}\,\, \mathbb{R}_{+}^{n+1}\, \Big(\mathbb{R}_{-}^{n+1}\Big)  \Longleftarrow\!\!\Longrightarrow \Phi^{*}\,\, \mathrm{\it
is \,\,(left)\,\,regular\,\,on}\,\, \mathbb{R}_{-}^{n+1}\,\Big(\mathbb{R}_{+}^{n+1}\Big).
\eeq
Hence the symmetry principle for regular functions is proved.
\ef

By Remark \ref{re 4.1}, we may get the following  results for boundary values.
\ble\label{le 4.3}
If $\Phi$ is defined in $\mathbb{R}_{+}^{n+1}$ with the boundary value $\Phi^{+}(t)$, then $\Phi^{*}$  has the boundary value $\big[\Phi^{*}\big]^{-}(t)$  and
\beq\label{4.44}
\big[\Phi^{+}\big]^{*}(t)=\big[\Phi^{*}\big]^{-}(t),
\,\,\,\,\,\,t\in\mathbb{R} _{0}^{n+1}.
\eeq
In other words,
\beq\label{4.45}
\Phi^{+}(t)+\left[\Phi^{*}\right]^{-}(t)=2\, \mathrm{Re}\!\left(\Phi^{+}(t)\right),\,\,\,\,t\in\mathbb{R} _{0}^{n+1}.
\eeq
\ele
\ble\label{le 4.4}
If $\Phi$ is defined in $\mathbb{R}_{-}^{n+1}$ with the boundary value $\Phi^{-}(t)$, then $\Phi^{*}$ also has the boundary value $\big[\Phi^{*}\big]^{+}(t)$  and
\beq\label{4.46}
\big[\Phi^{-}\big]^{*}(t)=\big[\Phi^{*}\big]^{+}(t),\,\,\,\,t\in\mathbb{R} _{0}^{n+1}.
\eeq
In other words,
\beq\label{4.47}
\Phi^{-}(t)+\left[\Phi^{*}\right]^{+}(t)=2\, \mathrm{Re}\!\left(\Phi^{-}(t)\right),\,\,\,\,t\in\mathbb{R} _{0}^{n+1}.
\eeq
\ele

By Lemma  \ref{th 4.1} and Lemma \ref{le 4.3}, we immediately have the following theorem, which is the basis for the Hilbert boundary value problem being transferable  into the Riemann boundary value problem.\vspace{1.5mm}

If $\Phi$ is defined in $\mathbb{R} _{+}^{n+1}$\, $\Big(\mathbb{R} _{-}^{n+1}\Big)$, we call the function
\beq\label{4.48}
\big(\mathcal{E}[\Phi]\big)(w)=\left\{
\begin{array}{ll}
\Phi(w), &{\it when}\,\,\,\,w\in\mathbb{R}_{+}^{n+1}\,\,\Big(\mathbb{R} _{-}^{n+1}\Big),\\[3.5mm]
{\Phi}^{*}(w),&{\it when}\,\,\,\,w\in\mathbb{R}_{-}^{n+1}\,\,
\Big(\mathbb{R} _{+}^{n+1}\Big)
\end{array}
\right.
\end{equation}
the
symmetric extension of $\Phi$ defined in $\mathbb{R}^{n+1}_{+}$ $\Big(\mathbb{R} _{-}^{n+1}\Big)$.

\bth\label{th 4.2}
If $\Phi$ is regular on $\mathbb{R}_{+}^{n+1}\,\Big(\mathbb{R} _{-}^{n+1}\Big)$
and can be extended \vspace{1mm}
continuously to $\mathbb{R}_{0}^{n+1}$, then its
symmetric extension $\mathcal{E}[\Phi]$ is the sectionally                                     {holomorphic function} with $\mathbb{R}_{0}^{n+1}$ {as the} jump surface
$($\mbox{see \rm\cite{ddq17}}$)$.
\eeth

With the help of this theorem, the Hilbert BVP (\ref{3.20}) will be converted equivalently to the Riemann BVP discussed in \cite{ddq17} with an additional restricting condition.
\bth\label{th 4.3}
The Schwarz BVP $(\ref{3.22})$ is equivalent to the following
Riemann BVP $(\ref{4.51})$  under the relationship
\beq\label{4.49}
\Psi(w)=\Big(\mathcal {E} [\Phi]\Big)(w),\,\,\,\,\,\,\,\,\,\,w\in\mathbb{R}_{+}^{n+1}\cup
\mathbb{R}_{-}^{n+1}
\eeq
or
\beq\label{4.50}
\Phi(w)=\Psi\Big|_{\mathbb{R}_{+}^{n+1}}(w)=\Psi^{+}(w),\,\,\,\,\,\,w\in\mathbb{R}_{+}^{n+1}.
\eeq
\eeth

\mbox{\bf $R_m^{*}$ problem with the reflection condition.} \vspace{1mm}Find a sectionally {holomorphic} function $\Psi$, with $\mathbb{R} ^{n+1}_{0}$ as its jump plane  such that \beq\label{4.51} \left\{
\begin{array}{l}
\Psi^{+}(x)+\Psi^{-}(x)=2c(x),\,\, x\in\mathbb{R} ^{n+1}_{0}\,\,\,\,(\mbox{\it boundary value condition}),\\[2.6mm]
\Psi(w)=o\big(w^{m+1}\big)\,\,\,\mbox{\it near}\,\,\,\infty\,\,\,
(\mbox{\it growth condition at the infinity} ),\\[3mm]
\big[\Psi^{+}\big]^{*}(x)=\big[\Psi^{*}\big]^{-}(x),\,\,\,\,\,\,x\in\mathbb{R} _{0}^{n+1}
\,\,\,\,(\mbox{\it reflection condition}),
\end{array}
\right.
\eeq
where $m$ is some integer.

\vspace{0.5mm}

\mbox{\bf The proof of Theorem \ref{th 4.3}.}\,\, Let $\Phi$ be a solution of the Schwarz problem (\ref{3.22}), by using Theorem \ref{th 4.1}, (\ref{4.20}) and (\ref{4.45}), we know that $\Psi(w)=\big(\mathcal {E} [\Phi]\big)(w)$ is the solution of the $R^{*}_m$ problem (\ref{4.51}).
And, in turn, if $\Psi$ is the solution {of the} $R^{*}_m$ problem (\ref{4.51}), then
\beq\label{4.52}
\Phi(w)=\Phi^{+}(w)=\Psi^{+}(w),\,\,\,\,w\in \mathbb{R} _{+}^{n+1}
\eeq
is the solution of the Schwarz problem (\ref{3.22}). In fact, by the reflection condition
in (\ref{4.51})  we have
(\ref{4.47}) which results in (\ref{3.22}). \ef

Now, the remaining question is how to
solve the $R_m^{*}$ problem  (\ref{4.51}). \vspace{1mm}
To do so, we introduce the self-reflex action of $\Phi$ defined $\mathbb{R}^{n+1}_{+}\bigcup\mathbb{R}^{n+1}_{-}$.\vspace{1mm}

\subsection{Self-reflex action}
\mbox{}\hspace{6mm}
If $\mho$ is defined on $\mathbb{R}_{+}^{n+1}\bigcup\mathbb{R}_{-}^{n+1}$, say,
\beq \label{4.53}
\mho(w)=\left\{\begin{array}{ll}
\mho^{+}(w), &w\in \mathbb{R}_{+}^{n+1},\\[2mm]
\mho^{-}(w),&w\in\mathbb{R}_{-}^{n+1},
\end{array}
\right.
\eeq
then
\beq \label{4.54}
\mho_{\updownarrow }(w)=\left\{\begin{array}{ll}
\big[\mho^{-}\big]^{*}(w), &w\in  \mathbb{R}_{+}^{n+1},\\[2.5mm]
\big[\mho^{+}\big]^{*}(w),&w\in \mathbb{R}_{-}^{n+1},
\end{array}
\right.
\eeq
is called the reflective function of $\mho$. Obviously, by  (\ref{4.17})
\beq\label{4.55}
\mho(w)=\big(\mho_{\updownarrow }\big)_{\updownarrow }(w),\,\,\,\,w\in\mathbb{R}_{\pm}^{n+1}.
\eeq
In particular, if
\beq\label{4.56}
\mho_{\updownarrow }(w)=\mho(w),\,\,\,\,w\in\mathbb{R}_{\pm}^{n+1},
\eeq
then we  call $\mho$ a self-reflection function. Obviously, by (\ref{4.17}),
\beq \label{4.57}
\big(\mathscr{R}[\mho]\big)(w)=\displaystyle\frac{\mho(w)+\mho_{\updownarrow }(w)}{2}
\eeq
is a self-reflection function, which is called the self-reflection function of $\mho$.
For the sake of convenience, we call the above steps from $\Phi,$ defined on $\mathbb{R}^{n+1}_{+}$ $\left(\mathbb{R}^{n+1}_{+}\right)$,
to $\mathscr{R}\big[\mathcal{E}[\Phi]\big]$  to be the so-called self-reflex action of $\Phi$.
From (\ref{4.56}) and (\ref{4.57}) we know that
$\Phi$ is a self-reflection function  if and only if
\beq\label{4.58}
\Phi(w)=\big(\mathscr{R}[\Phi]\big)(w),\,\,\,\, w\in\mathbb{R}_{+}^{n+1}\cup\mathbb{R}_{-}^{n+1}.
\eeq
\bex\label{ex 4.1}{\rm From Example \ref{ex 2.1}, we know that
 the restriction $E\big|_{\mathbb{R}^{n+1}_{\pm}}$ of the Cauchy kernel given in $(\ref{2.34})$ is a self-reflection
function.}
\eex
\bex\label{ex 4.2} {\rm Assume $f$ is a {para real-valued function}.
Let
\beq\label{4.59}
\Psi(w)=\left\{
\begin{array}{ll}
\Big(\mathcal {S}[f]\Big)(w)=
 \displaystyle\frac{\,1}{\bigvee_{n+1}}\!
 \int_{\mathbb{R}^{n+1}_{0}}E(x-w)\,\mathrm{d}\sigma f(x),\,\,\,\,\,
&w\!\in\!\mathbb{R} ^{n+1}_{+},\\[7mm]
-\Big(\mathcal {S}[f]\Big)(w)=
 -\displaystyle\frac{\,1}{\bigvee_{n+1}}\!
 \int_{\mathbb{R}^{n+1}_{0}}E(x-w)\,\mathrm{d}\sigma f(x),\,\,\,\,\,
&w\!\in\!\mathbb{R} ^{n+1}_{-},
\end{array}
\right.
\eeq
where $\mathcal {S}[f]$ is Cauchy type integral on the hyperplane
$\mathbb{R} _{0}^{n+1}$,
and
\beq\label{4.60}
\mathrm{d}S=\mathrm{d}x_0\mathrm{d}x_1\cdots\mathrm{d}x_{n-1}=e_n\,\mathrm{d}\sigma
\eeq
is the elementary surface measure on the hyperplane $\mathbb{R}^{n+1}_{0}$ $\big($see \cite{ddq17}$\big)$.
Then, by  Example \ref{ex 4.1} and Lemma \ref{4.1},
its self-reflection function is
\beq\label{4.61}
\Big(\mathscr{R}\big[\Psi\big]\Big)(w)=\Big(\mathcal {S}[f]\Big)(w)=
 \displaystyle\frac{\,1}{\bigvee_{n+1}}\!
 \int_{\mathbb{R}^{n+1}_{0}}E(x-w)\,\mathrm{d}\sigma f(x),\,\,\,\,\,\,\,\,
w\in\mathbb{R} _{\pm}^{n+1}.
\eeq
}
\eex

\bex\label{ex 4.3}{
\rm Let
\beq\label{4.62}
\Psi(w)=\left\{
\begin{array}{ll}
\,\,Z^{\alpha}(w)\,c_{\alpha},\,\,\,\,\,& w\in\mathbb{R}_{+}^{n+1},\\[2mm]
\!\!\!-Z^{\alpha}(w)\,c_{\alpha}\,\,\,\,\, &w\in\mathbb{R}_{-}^{n+1},
\end{array}
\right.
\eeq
where $c_{\alpha}$ is a hypercomplex constant and $Z^{\alpha}$ is a Fueter polynomail given in (\ref{3.10}). Then,
\beq\label{4.63}
\Big(\mathscr{R}\big[\Psi\big]\Big)(w)=Z^{\alpha}(w)\, \mathrm{Im\,}c_{\alpha}\, e_n,
\,\,\,\,\,\,\,\,
w\in\mathbb{R} _{\pm}^{n+1},
\eeq
which is called the para-imaginary coefficient polynomial.
}
\eex

We point out an obvious fact that,
 if $\Psi$ is the solution of $R_m^{*}$  problem (\ref{4.51}) with a reflection condition,
  of course it is the solution of
 the following Riemann boundary value problem $R_m$ discussed in \cite{ddq17}.
 \vspace{1mm}

 \mbox{\bf $R_m$ problem with no reflection condition.}  Find a  sectionally {holomorphic}  function $\Psi$, with $\mathbb{R} ^{n+1}_{0}$ as its jump plane,  such that
 \beq\label{4.64} \left\{
\begin{array}{l}
\Psi^{+}(x)+\Psi^{-}(x)=2c(x),\,\, x\in\mathbb{R} ^{n+1}_{0}\,\,\,\,(\mbox{\it boundary value condition}),\\[2.6mm]
\Psi(w)=o\big(w^{m+1}\big)\,\,\,\mbox{\it near}\,\,\,\infty\,\,\,
(\mbox{\it growth condition at the infinity}).
\end{array}
\right.
\eeq

If $\Psi$ is both a solution of (\ref{4.64}) and a reflexive function, then it is called a reflexive solution of (\ref{4.64}). Under such case, by (\ref{4.55})
the reflection condition in (\ref{4.51}) automatically holds.
So, we  have the following result.
\ble\label{le 4.5}
The reflexive solution {of the} $R_m$ problem  with no reflection condition  $(\ref{4.64})$
{surely} is the solution {of the} $R^{*}_m$ problem
 with the reflection condition $(\ref{4.51})$.
\ele

  In {turn}, we also have the following result.
\ble\label{le 4.6} If $\Psi$ is the solution {of the} $R_m$ problem without reflective condition $(\ref{4.64})$,  then its
reflection function $\Psi_{\updownarrow}$ and
the
self-reflection function $\mathscr{R}[\Psi]$ are
  the solutions of $(\ref{4.51})$.
\ele

Thus, $\mathscr{R}[\Psi]$ is  the self-reflection solution of  the  $R^{*}_m$ problem with reflection condition  in $(\ref{4.51})$.
\vspace{0.1mm}

\p\,\,First, the reflection function $\Psi_{\updownarrow}$ is regular on $\mathbb{R}_{\pm}^{n+1}$ by Lemma \ref{le 4.2}.  Secondly,
using the reflection operator to both sides of the formulas in (\ref{4.64}) we have, by
(\ref{4.20}) and (\ref{4.9}),
 \beq\label{4.65}
 \left\{
\begin{array}{l}
\left(\Psi_{\updownarrow}\right)^{+}(x)+\left(\Psi_{\updownarrow}\right)^{-}(x)=2c(x),\,\, x\in\mathbb{R} ^{n+1}_{0}\,\,\,\,(\mbox{\it boundary value condition}),\\[2.6mm]
\Psi_{\updownarrow}(w)=o\big(w^{m+1}\big)\,\,\,\mbox{\it near}\,\,\,\infty\,\,\,
(\mbox{\it growth condition at the infinity}).
\end{array}
\right.
\eeq
This is to say that $\Psi_{\updownarrow}$ is also a solution of (\ref{4.64}). So is the reflective function
$
 \mathscr{R}[\Psi]$  given in (\ref{4.57}) with $\mho=\Psi$,
 since $\Psi$ and $\Psi_{\updownarrow}$ are the solutions of (\ref{4.64}). \ef

\bth
\label{th 4.4}
 The general solution of the  Schwarz boundary value problem  $(\ref{3.22})$
should be
 \begin{equation}\label{4.66}
 \Phi(w)=\Big(\mathscr{R}[\Psi]\Big)(w)=\displaystyle\frac{\Psi(w)+{\Psi}_{\updownarrow }(w)}{2},\,\,\,\,\,\,\,w\in \mathbb{R}^{n+1}_{+},
\end{equation}
where $\Psi$ is the solution of the ${R}_{m}$ problem $(\ref{4.64})$.
\eeth

Thus, we have the principle of so-called self-reflex action  by Theorem \ref{th 4.2},   Lemma \ref{le 4.5} and Lemma \ref{le 4.6}.  In short, the solutions of
$R_m$ problem $(\ref{4.46})$  with no reflection condition
are  derived from the solutions of the Schwarz boundary value problem  $(\ref{3.22})$
 by the self-reflex action, or the solutions of the Schwarz boundary value problem  $(\ref{3.22})$ may be obtained from  the  solutions of $R_m$ problem  without reflection condition  $(\ref{4.64})$ through two steps: \vspace{1mm} first taking the self-reflex action and then taking the restriction on $\mathbb{R}_{+}^{n+1}$.

\mysection{Solutions of $\mathrm{H}_m$ problem}

Based on the reflexive principle, in order to solve $S_m$ problem $(\ref{3.22})$ we only need to solve the $R_m$ problem $(\ref{4.64})$, which is discussed in detail in\cite{ddq17}. \vspace{2.5mm}

\subsection{Solutions of $R_m$ promlem}
\mbox{}\hspace{6mm}For the convenience of reference, here we restate the results for $R_m$ problem as follows.
\bth [\mbox{\rm see \cite{ddq17}}]\label{th 5.1}
For the Riemann boundary value problem $R_{m}$ $(\ref{4.64})$ the following four cases are a complete classification.
\vspace{1.5mm}

$(1)$ Let $m\geq 0$, $c\in \widehat{H}^{\mu}\Big(\mathbb{R}_{0}^{n+1}\Big)$, then its general
solution is
\beq\label{5.1}
\Psi(w)=\left\{
\begin{array}{ll}
\,\,\Phi(w),&w\in\mathbb{R}_{+}^{n+1},\\[2mm]
\!\!\!-\Phi(w),&w\in\mathbb{R}_{-}^{n+1},
\end{array}\right.
\eeq
where
\beq
\label{5.2}
\Phi(w)\!=\!\Big(\mathcal {S}[c]\Big)(w)\!+\!P_{m}(w)\!=\!
 \displaystyle\frac{1}{\bigvee_{n+1}}\!\int_{\mathbb{R}^{n+1}_{0}}E(x-w)\mathrm{d}\sigma c(x)
  \!+\!\displaystyle\sum_{|\alpha|=0}^{m}\!\displaystyle\frac{1}{|\alpha|!}\, Z^{\alpha}(w)\,c_{\alpha},\,\,
w\!\in\!\mathbb{R} ^{n+1}_{\pm},
\eeq
where  $P_m$ is  arbitrary hypercomplex symmetric polynomial of degree not exceeding $m$
with
$C^{m}_{n+m}$ free hypercomplex constants $c_{\alpha}$.
\vspace{2.5mm}

$(2)$ Let $m=-1$ with\vspace{2mm} $c\in\widehat{H}\left(\mathbb{R}^{n+1}_{0}\right)$,
it has
the unique solution
\beq\label{5.3}
\Phi(w)=\left\{
\begin{array}{ll}
\,\Big(\mathcal {S}[c]\Big)(w)=
 \displaystyle\frac{1}{\bigvee_{n+1}}\int_{\mathbb{R}^{n+1}_{0}}E(x-w)
 \,\mathrm{d}\sigma\, c(x),&
w\in\mathbb{R} ^{n+1}_{+},\\[8mm]
\!\!\!-\Big(\mathcal {S}[c]\Big)(w)=-\displaystyle\frac{1}{\bigvee_{n+1}}\int_{\mathbb{R}^{n+1}_{0}}E(x-w)
 \,\mathrm{d}\sigma\, c(x),&
w\in\mathbb{R} ^{n+1}_{-},
\end{array}
\right.
\eeq
if and only if
\beq\label{5.4}
c(\infty)=:\lim_{x\in \mathbb{R} ^{n+1}_{0},\,\, x\rightarrow\infty}\,c(x)=0.
\eeq
\vspace{1.5mm}

$(3)$ \vspace{2mm}Let $-n\!<\!m\!<\!-1$ and  $r\!=\!-(m\!+\!1)$, with $c\!\in\! \widehat{H}_{r,0}\Big(\mathbb{R} ^{n+1}_{0}\Big)$, it has the unique solution
$(\ref{5.3})$.\vspace{2mm}

$(4)$ \vspace{2.6mm}
Let $m\leq-n$ and  $r\!=\!-(m\!+\!1)$, with $c\in \widehat{H}_{r,0}\Big(\mathbb{R} ^{n+1}_{0}\Big)$,  it has the unique solution $(\ref{5.3})$
if the $C^{n}_{-m-1}$ conditions
\beq\label{5.5}
          \int_{\mathbb{R}_{0}^{n+1}}Z^{\alpha}(x)\,\mbox{\rm d}\sigma \,c(x)=0,\,\,\,|\alpha|=0,1,\cdots\!, -(n+1+m)
\eeq
are fulfilled.
\eeth

\bre\label{re 5.1}
The $R_m$ problem discussed in\mbox{\rm \cite{ddq17}}
 is the jump problem for $\Phi$. Here the $R_m$ problem $(\ref{4.64})$ for $\Psi$ is called the \mbox{\rm Szeg\"{o}} problem in some literature. They are slightly different and governed  by the relation $(\ref{5.1})$.
\ere

 \subsection{Solutions of $\mathrm{S}_m$ promlem}
\mbox{}\hspace{6mm}\vspace{1mm}
By Theorem \ref{4.4}, Theorem \ref{5.1}, Example \ref{ex 4.2} and
Example \ref{ex 4.3},   we get the following result.
\bth\label{th 5.2}
The general solution of the  Schwarz boundary value problem  $(\ref{3.22})$
should be as following four cases.
\vspace{2mm}

$\mathbf{Case\, 1.}$ When $m\geq 0$ and $c\in \widehat{H}^{\mu}\Big(\mathbb{R}_{0}^{n+1}\Big)$, then
\beq\label{5.6}
\Phi(w)=\Big(\mathcal {S}[c]\Big)(w)=
 \displaystyle\frac{\,1}{\bigvee_{n+1}}\!\int_{\mathbb{R}^{n+1}_{0}}
 E(x-w)\mathrm{d}\sigma c(x)+\displaystyle\sum_{|\alpha|=0}^{m}\!\displaystyle\frac{1}{|\alpha|!}\, Z^{\alpha}(w)\,R_{\alpha}\,e_n\,\,\,\,\,\,\,
w\in\mathbb{R}_{+}^{n+1},
\eeq
where $R_{\alpha}$ are $C^{m}_{n+m}$ free para-real hypercomplex constants.
\vspace{3mm}

$\mathbf{Case\, 2.}$ When $m=-1$ with\vspace{2mm}  $c\in\widehat{H}\Big(\mathbb{R}^{n+1}_{0}\Big)$,
it has the unique solution
\beq\label{5.7}
\Phi(w)=\Big(\mathcal {S}[c]\Big)(w)=
 \displaystyle\frac{1}{\bigvee_{n+1}}\int_{\mathbb{R}^{n+1}_{0}}E(x-w)
 \,\mathrm{d}\sigma\, c(x),\,\,\,\,
w\in\mathbb{R} ^{n+1}_{+},
\eeq
if and only if
\beq\label{5.8}
c(\infty)=\lim_{x\in \mathbb{R} ^{n+1}_{0},\,\, x\rightarrow\infty}c(x)=0.
\eeq
\vspace{3mm}

$\mathbf{Case\, 3.}$  When $-n\!<\!m\!<\!-1$ and  $r\!=\!-(m\!+\!1)$,
with $c\!\in\! \widehat{H}_{r,0}\Big(\mathbb{R} ^{n+1}_{0}\Big)$, it has the unique solution
$(\ref{5.7})$.\vspace{3mm}

$\mathbf{Case\, 4.}$ \vspace{2.6mm}
When $m\leq-n$ and  $r\!=\!-(m\!+\!1)$, with $c\in \widehat{H}_{r,0}\Big(\mathbb{R} ^{n+1}_{0}\Big)$,
it has the
unique solution $(\ref{5.7})$
if the $C^{n}_{-m-1}$ conditions
\beq\label{5.9}
                   \int_{\mathbb{R}_{0}^{n+1}}Z^{\alpha}(x)\,\mbox{\rm d}\sigma\, c(x)=0,\,\,\,|\alpha|=0,1,\cdots\!, -(n+1+m)
\eeq
or

\beq\label{5.10}
          \int_{\mathbb{R}_{0}^{n+1}}Z^{\alpha}(x)\, c(x)\,\mbox{\rm d}x_1 \mbox{\rm d}x_2\cdots \mbox{\rm d}x_n=0,\,\,\,|\alpha|=0,1,\cdots\!, -(n+1+m)
\eeq
are fulfilled.

\eeth

\subsection{Solutions of $\mathrm{H}_m$ promlem}
\mbox{}\hspace{6mm}
The Hilbert problem (\ref{3.20}) for the function $\Phi$ may be directly translated into the Schwarz problem for the function
$\Psi$ by
\beq\label{5.11}
\Phi(w)=\left\{
\begin{array}{ll}
\,\Psi(w),&w\in\mathbb{R}_{+}^{n+1},\\[2mm]
\Psi(w)\lambda,&w\in\mathbb{R}_{-}^{n+1}.
\end{array}
\right.
\eeq
\bth\label{th 5.3}
The general solution of the  Hilbert boundary value problem  $(\ref{3.20})$
should be as following four cases. \vspace{3mm}

$\mathbf{Case\, 1.}$ When $m\geq 0$ and $c\in \widehat{H}^{\mu}\Big(\mathbb{R}_{0}^{n+1}\Big)$, then
\beq\label{5.12}
\Phi(w)=
 \displaystyle\frac{\,1}{\bigvee_{n+1}}\!\int_{\mathbb{R}^{n+1}_{0}}E(x-w)\mathrm{d}
 \sigma c(x)\lambda^{-1}+\displaystyle\sum_{|\alpha|=0}^{m}\!\displaystyle\frac{1}{|\alpha|!}\, Z^{\alpha}(w)\,R_{\alpha}\,\lambda^{-1},
w\in\mathbb{R}_{+}^{n+1},
\eeq
where $R_{\alpha}$ are $C^{m}_{n+m}$ free para-real hypercomplex constants.
\vspace{3mm}

$\mathbf{Case\, 2.}$ When $m=-1$ with\vspace{2mm} $c\in\widehat{H}\left(\mathbb{R}^{n+1}_{0}\right)$,
it has the unique solution
\beq\label{5.13}
\Phi(w)=\Big(\mathcal {S}[c]\Big)(w)\lambda^{-1}=
 \displaystyle\frac{1}{\bigvee_{n+1}}\int_{\mathbb{R}^{n+1}_{0}}E(x-w)
 \,\mathrm{d}\sigma\, c(x)\, \lambda^{-1},\,\,\,\,
w\in\mathbb{R} ^{n+1}_{+},
\eeq
if and only if
\beq\label{5.14}
c(\infty)=\lim_{x\in \mathbb{R} ^{n+1}_{0},\,\, x\rightarrow\infty}c\!(x)=0.
\eeq
\vspace{1mm}

$\mathbf{Case\, 3.}$  When $-n\!<\!m\!<\!-1$ and  $r\!=\!-(m\!+\!1)$, with  $c\!\in\! \widehat{H}_{r,0}\Big(\mathbb{R} ^{n+1}_{0}\Big)$, it has the unique solution
$(\ref{5.13})$.
\vspace{3mm}

$\mathbf{Case\, 4.}$ \vspace{2.6mm}
Let $m\leq-n$ and  $r\!=\!-(m\!+\!1)$. If $c\in \widehat{H}_{r,0}\Big(\mathbb{R} ^{n+1}_{0}\Big)$, then it has the
unique solution $(\ref{5.13})$
provided the $C^{n}_{-m-1}$ conditions
\beq\label{5.15}
                   \int_{\mathbb{R}_{0}^{n+1}}Z^{\alpha}(x)\,\mbox{\rm d}\sigma\, c(x)=0,\,\,\,|\alpha|=0,1,\cdots\!, -(n+1+m)
\eeq
or
\beq\label{5.16}
          \int_{\mathbb{R}_{0}^{n+1}}Z^{\alpha}(x)\,c(x)\,\mbox{\rm d}x_1 \mbox{\rm d}x_2\cdots \mbox{\rm d}x_n=0,\,\,\,|\alpha|=0,1,\cdots\!, -(n+1+m)
\eeq
are fulfilled.

\eeth
\vspace{4mm}

%



\vspace{32mm}
\end{document}